\documentclass[12pt]{amsart}

 \usepackage{amsfonts,graphics,amsmath,amsthm,amsfonts,amscd,amssymb,amsmath,latexsym,multicol}
\usepackage{epsfig,url}
\usepackage{flafter}

\makeatletter

\def\jobis#1{FF\fi
  \def\predicate{#1}%
  \edef\predicate{\expandafter\strip@prefix\meaning\predicate}%
  \edef\job{\jobname}%
  \ifx\job\predicate
}

\makeatother

\if\jobis{proposal}%
\else
\fi

 \usepackage[matrix, arrow]{xy}

\DeclareMathOperator{\Supp}{Supp}

 \newcommand{\N}{\mathbb N}
 
 \newcommand{\Q}{\mathbb Q}
 \newcommand{\R}{\mathbb R}
 \newcommand{\Z}{\mathbb Z}
 \newcommand{\bir}{\dashrightarrow}
 \newcommand{\rddown}[1]{\left\lfloor{#1}\right\rfloor} 
  \newcommand{\ep}{\varepsilon}

 \numberwithin{equation}{subsection}
 \numberwithin{footnote}{subsection}

 \newtheorem{cor}[subsection]{Corollary}
 \newtheorem{lem}[subsection]{Lemma}
 \newtheorem{prop}[subsection]{Proposition}
 \newtheorem{thm}[subsection]{Theorem}
 \newtheorem{conj}[subsection]{Conjecture}

{
    \newtheoremstyle{upright}%
        {8pt plus2pt minus4pt}%
        {8pt plus2pt minus4pt}%
        {\upshape}%
        {}%
        {\bfseries\scshape}%
        {}%
        {1em}%
        {}%

\theoremstyle{upright}

 \newtheorem{defn}[subsection]{Definition}

 \newtheorem{rem}[subsection]{Remark}

}

 \newcommand{\ke}[1]{$\acute{\mbox{e}}$}
 \newcommand{\ku}[1]{$\acute{\mbox{u}$}}
 \newcommand{\kl}[1]{$\acute{\mbox{l}}$}
 \newcommand{\kh}[1]{$\acute{\mbox{h}}$}
 \newcommand{\kr}[1]{$\acute{\mbox{r}}$}
 \newcommand{\kx}[1]{$\acute{\mbox{x}}$}
 \newcommand{\ki}[1]{${\^\i}$}

\title{On existence of log minimal models II}
\author{Caucher Birkar}
\date{\today}
\begin{document}
\maketitle

\begin{abstract}
We prove that the existence of log minimal models in dimension $d$ essentially implies 
the LMMP with scaling in dimension $d$. As a consequence we prove that a weak
nonvanishing conjecture in dimension $d$ implies the minimal model conjecture
in dimension $d$.
\end{abstract}


\section{Introduction}

We work over a fixed algebraically closed field $k$ of characteristic zero. See section $2$ for notation and terminology.
Remember that a lc pair $(X/Z,B)$ is called pseudo-effective if $K_X+B$ is pseudo-effective$/Z$, that is, 
if there is a sequence of $\R$-divisors $M_i\ge 0$ such that $K_X+B\equiv \lim_{i\to \infty} M_i$ in $N^1(X/Z)$. 
The pair is called effective if $K_X+B\equiv M/Z$ for some $M\ge 0$.

The following two conjectures are, at the moment, the most important open problems in birational geometry and the classification theory of 
algebraic varieties. 

\begin{conj}[Minimal model]\label{mmodel}
Let $(X/Z,B)$ be a lc pair. If it is pseudo-effective then it has a log minimal model, and if it is not pseudo-effective 
then it has a Mori fibre space.
\end{conj}

\begin{conj}[Abundance]\label{abundance}
Let $(X/Z,B)$ be a lc pair. If $K_X+B$ is nef$/Z$, then it is semi-ample$/Z$.
\end{conj}

For a brief history of the many results on the minimal model conjecture see the introduction to [\ref{B}]. 
On the other hand, there has been little progress regarding the abundance conjecture in higher dimension. 
The main conceptual obstacle to abundance is the following problem.

\begin{conj}[Weak nonvanishing]\label{c-wa}
Let $(X/Z,B)$ be a $\Q$-factorial dlt pair. If $K_X+B$ is pseudo-effective$/Z$, then it is effective$/Z$, that is, 
$K_X+B\equiv M/Z$ for some $M\ge 0$.
\end{conj}

This conjecture is not only at the heart of the abundance conjecture but it is also closely related to the 
minimal model conjecture. In fact, we show that it implies the minimal model conjecture. 

\begin{thm}\label{t-main}
Assume the weak nonvanishing conjecture (\ref{c-wa}) in dimension $d$. Then,  the minimal model conjecture (\ref{mmodel}) 
holds in dimension $d$; moreover, if 
$(X/Z,B)$ is a $\Q$-factorial dlt pair of dimension $d$, then there is a sequence of divisorial 
contractions and log flips starting with $(X/Z,B)$ and ending up with a log minimal model or a Mori fibre 
space of $(X/Z,B)$. 
\end{thm}

The proof of this theorem is given via the following results and [\ref{B}, Proposition 3.4].

\begin{thm}\label{t-main2}
Assume the minimal model conjecture (\ref{mmodel}) in dimension $d$ for pseudo-effective $\Q$-factorial dlt pairs. 
Let $(X/Z,B+C)$ be a $\Q$-factorial lc pair of dimension $d$ such that 
\begin{enumerate}
\item $K_X+B+C$ is nef$/Z$,
\item $B,C\ge 0$, and
\item $(X/Z,B)$ is dlt.
\end{enumerate} 
Then, we can run the LMMP$/Z$ on $K_X+B$ with scaling of $C$, and it terminates if either 
\begin{itemize}
\item $B\ge H\ge 0$ for some ample$/Z$ $\R$-divisor $H$, or 
\item $C\ge H\ge 0$ for some ample$/Z$ $\R$-divisor $H$, or 
\item  $\lambda\neq \lambda_i$ for any $i$ where $\lambda$ and $\lambda_i$ 
are as in Definition \ref{d-scaling}.
\end{itemize}
\end{thm}

\begin{cor}\label{c-1}
Assume the minimal model conjecture (\ref{mmodel}) in dimension $d$ for pseudo-effective $\Q$-factorial dlt pairs. 
Let $(X/Z,B)$ be a $\Q$-factorial dlt pair of dimension $d$. Then, there is a sequence of divisorial 
contractions and log flips starting with $(X/Z,B)$ and ending up with a log minimal model or a Mori fibre 
space $(Y/Z,B_Y)$. In particular, the corresponding birational map $Y\bir X/Z$ 
does not contract divisors. 
\end{cor}

\begin{cor}\label{c-2}
Assume the minimal model conjecture (\ref{mmodel}) in dimension $d$ for pseudo-effective $\Q$-factorial dlt pairs.  
Then, the minimal model conjecture (\ref{mmodel}) holds in dimension $d+1$ for effective lc pairs. 
\end{cor}

We sometimes refer to some of the results of [\ref{BCHM}][\ref{BP}]. Actually, to prove the main results 
of this paper we only need two pages of 
[\ref{BP}], that is [\ref{BP}, Theorem 2.6], the rest that we need can be 
easily incorporated into the framework of [\ref{B}] and this paper. 


\section{Basics}

Let $k$ be an algebraically closed field of characteristic zero fixed throughout the paper. 

A \emph{pair} $(X/Z,B)$ consists of normal quasi-projective varieties $X,Z$ over $k$, an $\R$- divisor $B$ on $X$ with
coefficients in $[0,1]$ such that $K_X+B$ is $\mathbb{R}$-Cartier, and a projective 
morphism $X\to Z$. For a prime divisor $D$ on some birational model of $X$ with a
nonempty centre on $X$, $a(D,X,B)$
denotes the log discrepancy.

A pair $(X/Z,B)$ is called \emph{pseudo-effective} if $K_X+B$ is pseudo-effective/$Z$, that is,
up to numerical equivalence/$Z$ it is the limit of effective $\R$-divisors. The pair is called \emph{effective} if
$K_X+B$ is effective/$Z$, that is, there is an $\R$-divisor $M\ge 0$ such that $K_X+B\equiv M/Z$. 

 By a \emph{log flip}$/Z$ we mean the flip of a $K_X+B$-negative extremal flipping contraction$/Z$ for some 
lc pair $(X/Z,B)$ (cf. [\ref{B}, Definition 2.3]), 
and by a \emph{pl flip}$/Z$ we mean a log flip$/Z$ such that $(X/Z,B)$ is $\Q$-factorial dlt and the log flip is also an $S$-flip for 
some component $S$ of $\rddown{B}$.

A \emph{sequence of log flips$/Z$ starting with} $(X/Z,B)$ is a sequence $X_i\bir X_{i+1}/Z_i$ in which  
$X_i\to Z_i \leftarrow X_{i+1}$ is a $K_{X_i}+B_i$-flip$/Z$, $B_i$ is the birational transform 
of $B_1$ on $X_1$, and $(X_1/Z,B_1)=(X/Z,B)$.

In this paper, \emph{special termination} means termination near $\rddown{B}$ of any sequence of log flips$/Z$ 
starting with a pair $(X/Z,B)$, that is, 
the log flips do not intersect $\rddown{B}$ after finitely many of them. 

\begin{defn}\label{d-model}
A pair $(Y/Z,B_Y)$ is a \emph{log birational model} of $(X/Z,B)$ if we are given a birational map
$\phi\colon X\bir Y/Z$ and $B_Y=B^\sim+E$ where $B^\sim$ is the birational transform of $B$ and 
$E$ is the reduced exceptional divisor of $\phi^{-1}$, that is, $E=\sum E_j$ where $E_j$ are the
exceptional/$X$ prime divisors on $Y$. A log birational model $(Y/Z,B_Y)$ is a \emph{nef model} of 
$(X/Z,B)$ if in addition\\\\
(1) $(Y/Z,B_Y)$ is $\Q$-factorial dlt, and\\
(2) $K_Y+B_Y$ is nef/$Z$.\\

And  we call a nef model $(Y/Z,B_Y)$ a \emph{log minimal model} of $(X/Z,B)$ if in addition\\\\
(3) for any prime divisor $D$ on $X$ which is exceptional/$Y$, we have
$$
a(D,X,B)<a(D,Y,B_Y)
$$
\end{defn}

\begin{defn}[Mori fibre space]
A log birational model $(Y/Z,B_Y)$ of a lc pair $(X/Z,B)$ is called a Mori fibre space if 
$(Y/Z,B_Y)$ is $\Q$-factorial dlt, there is a $K_Y+B_Y$-negative extremal contraction $Y\to T/Z$ 
with $\dim Y>\dim T$, and 
$$
a(D,X,B)\le a(D,Y,B_Y)
$$
for  any prime divisor $D$ (on birational models of $X$) and the strict inequality holds if $D$ is on $X$ and contracted$/Y$.
\end{defn}

Our definitions of log minimal models and Mori fibre spaces are slightly different 
from the traditional ones, the difference being that we do not assume that $\phi^{-1}$ does not contract divisors. 
Even though we allow $\phi^{-1}$ to have exceptional divisors but these divisors are very special; if $D$ 
is any such prime divisor, then $a(D,X,B)=a(D,Y,B_Y)=0$.
Actually, in the plt case, our definition of log minimal models and the traditional one coincide (see [\ref{B}, Remark 2.6]).

\begin{defn}[LMMP with scaling]\label{d-scaling}
Let $(X_1/Z,B_1+C_1)$ be a lc pair such that $K_{X_1}+B_1+C_1$ is nef/$Z$, $B_1\ge 0$, and $C_1\ge 0$ is $\R$-Cartier. 
Suppose that either $K_{X_1}+B_1$ is nef/$Z$ or there is an extremal ray $R_1/Z$ such
that $(K_{X_1}+B_1)\cdot R_1<0$ and $(K_{X_1}+B_1+\lambda_1 C_1)\cdot R_1=0$ where
$$
\lambda_1:=\inf \{t\ge 0~|~K_{X_1}+B_1+tC_1~~\mbox{is nef/$Z$}\}
$$
 When $(X_1/Z,B_1)$ is $\Q$-factorial dlt, the last sentence follows from [\ref{B}, 3.1]. 
If $R_1$ defines a Mori fibre structure, we stop. Otherwise assume that $R_1$ gives a divisorial 
contraction or a log flip $X_1\bir X_2$. We can now consider $(X_2/Z,B_2+\lambda_1 C_2)$  where $B_2+\lambda_1 C_2$ is 
the birational transform 
of $B_1+\lambda_1 C_1$ and continue. That is, suppose that either $K_{X_2}+B_2$ is nef/$Z$ or 
there is an extremal ray $R_2/Z$ such
that $(K_{X_2}+B_2)\cdot R_2<0$ and $(K_{X_2}+B_2+\lambda_2 C_2)\cdot R_2=0$ where
$$
\lambda_2:=\inf \{t\ge 0~|~K_{X_2}+B_2+tC_2~~\mbox{is nef/$Z$}\}
$$
 By continuing this process, we obtain a sequence of numbers $\lambda_i$ and a 
special kind of LMMP$/Z$ which is called the \emph{LMMP$/Z$ on $K_{X_1}+B_1$ with scaling of $C_1$}; note that it is not unique. 
This kind of LMMP was first used by Shokurov [\ref{log-flips}].
When we refer to \emph{termination with scaling} we mean termination of such an LMMP. We usually put 
$\lambda=\lim \lambda_i$.

\emph{Special termination  with scaling} means termination near $\rddown{B_1}$ of any sequence of log flips$/Z$ with scaling 
of $C_1$, i.e. after finitely many steps, the locus of the extremal rays in the process do not intersect $\rddown{B_1}$.

When we have a lc pair $(X/Z,B)$, we can always find an ample$/Z$ $\R$-Cartier divisor $C\ge 0$ such that 
$K_X+B+C$ is lc and nef$/Z$,  so we can run the LMMP$/Z$ with scaling assuming that all the 
necessary ingredients exist, eg extremal rays, log flips. 
\end{defn}


\section{Extremal rays}

We need a result of Shokurov on extremal rays [\ref{ordered}]. Since we need stronger statements 
than those stated in [\ref{ordered}], we give 
detailed proofs here (see also [\ref{B2}]). Some parts of our proof are quite different from the originals. As a corollary, 
we give a short proof of a result of Kawamata on flops connecting minimal models.

Let $X\to Z$ be a projective morphism of normal quasi-projective varieties. A 
curve $\Gamma$ on $X$ is called \emph{extremal}$/Z$ if it generates an 
extremal ray $R/Z$ which defines a contraction $X\to S/Z$ and if for some ample$/Z$ divisor $H$ we 
have $H\cdot \Gamma=\min \{H\cdot \Sigma\}$ where $\Sigma$ ranges over curves generating $R$.
If $(X/Z,B)$ is dlt and $(K_X+B)\cdot R<0$, then by [\ref{antican}, Theorem] there is a curve $\Sigma$ generating $R$ 
such that $(K_X+B)\cdot \Sigma\ge -2\dim X$.
On the other hand, since $\Gamma$ and $\Sigma$ both generate $R$ we have 
$$
\frac{(K_X+B)\cdot \Gamma}{H\cdot \Gamma}=\frac{(K_X+B)\cdot \Sigma}{H\cdot \Sigma}
$$
hence 
\begin{equation}\label{eq-ac}
(K_X+B)\cdot \Gamma =((K_X+B)\cdot \Sigma) (\frac{H\cdot \Gamma}{H\cdot \Sigma})  \ge -2\dim X
\end{equation}\\

\begin{rem}\label{r-polytope}
Let $X/Z$ be a $\Q$-factorial dlt variety, $F$ be a reduced divisor on $X$, 
and ${V}$ be a rational affine subspace of the $\R$-vector space of divisors generated by the components of $F$. 
By [\ref{log-flips}, 1.3.2], the set
$$
\mathcal{L}=\{\Delta \in V \mid (X/Z,\Delta) ~~\mbox{is lc}\}
$$
is a rational polytope, that is, it is the convex hull of finitely many rational points in $V$. 
For any $\Delta\in \mathcal{L}$ and any extremal curve $\Gamma/Z$
the boundedness $(K_X+\Delta)\cdot \Gamma \ge -2\dim X$  holds as in \eqref{eq-ac}. Even though 
$(X/Z,\Delta)$ may not be dlt but we can use the fact that $(X/Z,a\Delta)$ is dlt for any 
$a\in [0,1)$.

Let $B_1,\dots, B_r$ be the vertices of $\mathcal{L}$, and let $m\in \N$ such that $m(K_X+B_j)$ are Cartier. 
For any $B\in \mathcal{L}$, there are  
 nonnegative real numbers $a_1,\dots, a_r$ such that $B=\sum a_jB_j$, $\sum a_j=1$, and each $(X/Z,B_j)$ is lc. 
Moreover, for any curve $\Gamma$ on $X$ 
the intersection number $(K_X+B)\cdot \Gamma$ can be written as $\sum a_j\frac{n_j}{m}$ for certain $n_1,\dots, n_r\in \Z$. 
If $\Gamma$ is extremal$/Z$, then the $n_j$ satisfy $n_j\ge -2m\dim X$.
\end{rem}

For an $\R$-divisor $D=\sum d_iD_i$ where the $D_i$ are the irreducible components of $D$, define $||D||:=\max\{|d_i|\}$.

\begin{prop}\label{p-rays}
Let $X/Z$, $F$, $V$, and $\mathcal{L}$ be as in Remark \ref{r-polytope}, and fix $B\in \mathcal{L}$. Then, 
there are real numbers $\alpha, \delta>0$, depending on $(X/Z,B)$ and $F$, such that 
\begin{enumerate}
\item if $\Gamma$ is any extremal curve$/Z$ and if 
$(K_X+B)\cdot \Gamma>0$, then $(K_X+B)\cdot \Gamma>\alpha$;

\item if $\Delta\in \mathcal{L}$, $||\Delta-B||<\delta$ and  $(K_X+\Delta)\cdot R\le 0$ for an extremal ray $R/Z$, 
then $(K_X+B)\cdot R\le 0$;

\item  let $\{R_t\}_{t\in T}$ be a family of extremal rays of $\overline{NE}(X/Z)$. Then, the set 
$$ 
\mathcal{N}_T=\{\Delta \in \mathcal{L} \mid (K_X+\Delta)\cdot R_t\ge 0 ~~\mbox{for any $t\in T$}\}
$$
is a rational polytope;

\item if $K_X+B$ is nef$/Z$, then for any $\Delta\in \mathcal{L}$ 
and for any sequence $X_i\bir X_{i+1}/Z_i$ of 
$K_X+\Delta$-flips$/Z$ which are flops with respect to $(X/Z,B)$ and any extremal curve $\Gamma/Z$ on $X_i$, if 
$(K_{X_i}+B_i)\cdot \Gamma>0$, then $(K_{X_i}+B_i)\cdot \Gamma>\alpha$ where $B_i$ is the 
birational transform of $B$;

\item assumptions as in $(5)$. In addition suppose that $||\Delta-B||<\delta$. If 
$(K_{X_i}+\Delta_i)\cdot R\le 0$ for an extremal ray $R/Z$ on some $X_i$, then $(K_{X_i}+B_i)\cdot R= 0$ 
where $\Delta_i$ is the birational transform of $\Delta$. 
\end{enumerate}
\end{prop}

{\textbf{Proof.}}
(1) If $B$ is a $\Q$-divisor, then the statement is trivially true even if $\Gamma$ is not extremal. If $B$ 
is not a $\Q$-divisor, 
let $B_1,\dots, B_r$, $a_1,\dots,a_r$, and $m$ be as in Remark \ref{r-polytope}. 
Then, 
$$
(K_X+B)\cdot \Gamma=\sum a_j(K_X+B_j)\cdot \Gamma
$$ 
and if $(K_X+B)\cdot \Gamma<1$, then there are only finitely many possibilities for the intersection 
numbers $(K_X+B_j)\cdot \Gamma$ because $(K_X+B_j)\cdot \Gamma \ge -2\dim X$. So, 
the existence of $\alpha$ is clear for (1).\\ 

(2)  If the statement is not 
true then there is an infinite sequence of $\Delta_t\in \mathcal{L}$ and extremal rays $R_t/Z$ such that for 
each $t$ we have 
$$
(K_X+\Delta_t)\cdot R_t\le 0  ~~~~~~~~ \mbox{,}~~~~ ~~~~~ (K_X+B) \cdot R_t> 0,
$$ 
and $||\Delta_t-B||$ converges to $0$. 
Let $B_1, \dots, B_r$ be the vertices of $\mathcal{L}$ which are 
rational divisors as $\mathcal{L}$ is a rational polytope. 
Then, there are nonnegative real numbers $a_1,\dots, a_r$ and $a_{1,t},\dots, a_{r,t}$
such that $B=\sum a_jB_j$, $\sum a_j=1$ and $\Delta_t=\sum a_{j,t} B_j$, $\sum a_{j,t}=1$. 
Since $||\Delta_t-B||$ converges to $0$, $a_j=\lim_{t\to \infty} a_{j,t}$. Perhaps after replacing the sequence with an infinite subsequence 
we can assume that the sign of $(K_X+B_j)\cdot R_t$ is independent of $t$, and that for each $t$ we have an extremal 
curve $\Gamma_t$ for $R_t$. Now, if $(K_X+B_j)\cdot \Gamma_t\le 0$, then  it 
is bounded from below hence there are only finitely many possibilities for this number and 
we could assume that it is independent of $t$. On the other hand, if $a_j\neq 0$, then  $(K_X+B_j)\cdot \Gamma_t$ 
 is bounded from below and above because 
$$
(K_X+\Delta_t)\cdot \Gamma_t=\sum a_{j,t}(K_X+B_j)\cdot \Gamma_t \le 0
$$ 
hence there are only finitely many possibilities for  $(K_X+B_j)\cdot \Gamma_t$  and 
we could assume that it is independent of $t$. 

Assume that $a_j\neq 0$ for $1\le j\le l$ but $a_j=0$ for 
$j>l$. Then, it is clear that 
$$
(K_X+\Delta_t)\cdot \Gamma_t=
$$
$$
 (K_X+B)\cdot \Gamma_t +\sum_{j\le l} (a_{j,t}-a_j) (K_X+B_j)\cdot \Gamma_t 
+\sum_{j>l} a_{j,t} (K_X+B_j)\cdot \Gamma_t
$$
would be positive by (1) if $t\gg 0$, which gives a contradiction.\\

(3) We may assume that for each $t\in T$ there is some $\Delta\in \mathcal{L}$ such that $(K_X+\Delta)\cdot R_t<0$, 
in particular, $(K_X+B_j)\cdot R_t<0$ for a vertex $B_j$ of $\mathcal{L}$. Since the set of such extremal rays is discrete, 
we may assume that $T\subseteq \N$.

Obviously, $\mathcal{N}_T$ is a convex compact subset of $\mathcal{L}$. If $T$ is finite, the claim is trivial. 
So we may assume that $T=\N$. By (2) and by the compactness of $\mathcal{N}_T$, there are 
$\Delta_1,\dots, \Delta_n\in \mathcal{N}_T$ 
and $\delta_1,\dots,\delta_n>0$ such that $\mathcal{N}_T$ is covered by 
$\mathcal{B}_i=\{\Delta \in \mathcal{L} \mid ||\Delta-\Delta_i||<\delta_i\}$ and such that 
if $\Delta\in \mathcal{B}_i$ with $(K_X+\Delta)\cdot R_t<0$ for some $t$, then $(K_X+\Delta_i)\cdot R_t=0$.
If
$$
T_i=\{t\in T \mid (K_X+\Delta)\cdot R_t<0 ~~\mbox{for some $\Delta \in \mathcal{B}_i$}\}
$$
then by construction $(K_X+\Delta_i)\cdot R_t=0$ for any $t\in T_i$.
Then, since the $\mathcal{B}_i$ give an open cover of $\mathcal{N}_T$, we have $\mathcal{N}_T=\bigcap_{1\le i\le n} \mathcal{N}_{T_i}$.
So, it is enough to prove that each $\mathcal{N}_{T_i}$ is a rational polytope and by replacing $T$ with $T_i$, 
we could assume from the beginning that there is some $\Delta\in \mathcal{N}_T$ such that 
$(K_X+\Delta)\cdot R_t=0$ for every $t\in T$.
 If $\dim \mathcal{L}=1$, this already proves the claim. 
If $\dim \mathcal{L}>1$, let $\mathcal{L}^1, \dots, \mathcal{L}^p$ be the proper faces of $\mathcal{L}$.
Then, each $\mathcal{N}_{T}^i=\mathcal{N}_{T}\cap \mathcal{L}^i$ is a rational polytope by induction. 
Moreover,  for each $\Delta''\in \mathcal{N}_{T}$ 
which is not $\Delta$, there is $\Delta'$ on some proper face of $\mathcal{L}$ such that $\Delta''$ is on the line 
segment determined by $\Delta$ and $\Delta'$. Since $(K_X+\Delta)\cdot R_t=0$ for every 
$t\in T$, if $\Delta'\in \mathcal{L}^i$, then $\Delta'\in \mathcal{N}_T^i$.
Hence $\mathcal{N}_T$ is the convex hull of $\Delta$ and all the  $\mathcal{N}_{T}^i$. 
Now, there is a finite subset $T'\subset T$ such that 
$$
\cup \mathcal{N}_T^i=\mathcal{N}_{T'}\cap (\cup \mathcal{L}^i)
$$
But then the convex hull of  $\Delta$ and $\cup \mathcal{N}_{T}^i$ 
is just $\mathcal{N}_{T'}$ and we are done.\\

(4) Since $K_X+B$ is nef$/Z$, $B\in \mathcal{N}_T$ where we take $\{R_t\}_{t\in T}$ to be the family of 
all the extremal rays of $\overline{NE}(X/Z)$. Since $\mathcal{N}_T$ is a rational polytope by (3), there are 
nonnegative real numbers $a_1',\dots,a_{r'}'$, and $m'\in \N$ so that $\sum a_j'=1$, $B=\sum a_j'B_j'$, 
and each $m'(K_X+B_j')$ is Cartier  
where $B_j'$ are the vertices of $\mathcal{N}_T$. Therefore, by the 
property $K_X+B=\sum a_j'(K_X+B_j')$, the sequence $X_i\bir X_{i+1}/Z_i$ 
is also a sequence of flops with respect to each $(X/Z,B_j')$. Moreover, $({X_i}/Z,B_{j,i}')$ is lc and $m'(K_{X_i}+B_{j,i}')$ is Cartier for 
any $j,i$ where $B_{j,i}'$ is the birational transform of $B_j'$. The rest is as in (1).\\

(5) Take $L$ to be the line in $V$ which goes through $B$ and $\Delta$ and let $\Delta'$ be the intersection 
point of $L$ and the boundary of $\mathcal{L}$, in the direction of $\Delta$. So, there are nonnegaitve real numbers 
$r,s$ such that $r+s=1$ and $\Delta=rB+s\Delta'$.  In particular, the sequence $X_i\bir X_{i+1}/Z_i$ is also a 
sequence of $K_X+\Delta'$-flips and $(X_i/Z,\Delta_i')$ is lc where $\Delta_i'$ is the birational transform of 
$\Delta'$. Suppose that there is an extremal ray $R/Z$ on some $X_i$ such that $(K_{X_i}+\Delta_i)\cdot R\le 0$ 
but $(K_{X_i}+B_i)\cdot R>0$.  Let $\Gamma$ be an extremal curve for $R$. By (4), $(K_{X_i}+B_i)\cdot \Gamma >\alpha$ 
and by \eqref{eq-ac} $(K_{X_i}+\Delta_i')\cdot \Gamma \ge -2\dim X$. Now 
$$
(K_{X_i}+\Delta_i)\cdot \Gamma=r(K_{X_i}+B_i)\cdot \Gamma+s(K_{X_i}+\Delta_i')\cdot \Gamma>r\alpha-2s\dim X
$$
and it is obvious that this is positive if $r>\frac{2s\dim X}{\alpha}$. In other words, if $\Delta$ is sufficiently close to 
$B$, then we get a contradiction. Therefore, it is enough to replace the $\delta$ of (2) by one sufficiently smaller. 
Note that we could also prove (2) in a similar way. 
$\Box$\\

In section 4, we will apply the proposition in a way similar to [\ref{ordered}].

Proposition \ref{p-rays} easily implies the following result of Kawamata [\ref{Kawamata}] on flops connecting log minimal models.

\begin{cor}\label{c-flops}
Let $(Y_1/Z,B_1)$ and $(Y_2/Z,B_2)$ be two klt pairs such that $K_{Y_1}+B_1$ 
and $K_{Y_2}+B_2$ are nef$/Z$, and $Y_1$ and $Y_2$  are isomorphic in codimension one. 
Then, $Y_1$ and $Y_2$ are connected by a sequence of flops$/Z$ with respect to $(Y_1/Z,B_1)$.
\end{cor}
{\textbf{Proof.}} 
Let $H_2$ be a general ample$/Z$ divisor on $Y_2$ and let $H_1$ be its birational transform on $Y_1$. 
There is $\delta>0$ such that $(Y_1/Z,B_1+\delta H_1)$ is klt. Now there is a general ample$/Z$ divisor 
$H_1'$ on $Y_1$ such that $(Y_2/Z,B_2+\delta H_2+\delta' H_2')$ is klt  for some $\delta'>0$ 
where $H_2'$ is the birational transform of $H_1'$. 
If $\delta$ is sufficiently small, then $K_{Y_1}+B_1+\delta H_1+\delta' H_1'$ is nef$/Z$. 
By [\ref{BCHM}][\ref{BP}], we can run the LMMP$/Z$ on $K_{Y_1}+B_1+\delta H_1$ with scaling of 
$\delta'H_1'$. After a finite sequence of log flips$/Z$, we end up with $Y_2$. On the other hand, we can lift the 
sequence to the $\Q$-factorial situation and by applying Proposition \ref{p-rays} we see that 
the sequence is a sequence of flops with respect to $(Y_1/Z,B_1)$ if $\delta$ is sufficiently small.
 $\Box$\\

Note that if  $(Y_1/Z,B_1)$ and $(Y_2/Z,B_2)$ are log minimal models of a klt pair $(X/Z,B)$, then 
$Y_1$ and $Y_2$ are automatically isomorphic in codimension one.


\section{Log minimal models and termination with scaling}

\vspace{0.5cm}
{\textbf{Proof of Theorem \ref{t-main2}}.}
\emph{Step 1.} The fact that we can run the LMMP$/Z$ on $K_X+B$ with scaling of $C$ follows from [\ref{B}, Lemma 3.1]. 
Note that the log flips required exist by the assumptions since existence of log flips is a special case 
of existence of log minimal models. Alternatively one can use  
[\ref{BCHM}][\ref{BP}].  We will deal with the termination statement. We may 
 assume that the sequence corresponding to the $\lambda_i$ is a sequence  $X_i\bir X_{i+1}/Z_i$ of log flips$/Z$  
starting with $(X/Z,B)$ where the $\lambda_i$
are obtained as in Definition \ref{d-scaling}. 
Remember that $\lambda=\lim_{i\to \infty} \lambda_i$.

If $B\ge H\ge 0$ for some ample$/Z$ $\R$-divisor 
$H$, then the LMMP terminates by  [\ref{BP}, Theorem 2.7]. Note that since $H$ 
is ample$/Z$, we can perturb the coefficients of $B$ and $C$ to reduce to the situation 
in which  $(X/Z,B+C)$ is klt 
(cf. [\ref{BP}, Remark 2.4]). If $C\ge H\ge 0$ where $H$ is an ample$/Z$ $\R$-divisor and if 
we have $\lambda>0$, then the termination 
follows again from [\ref{BP}, Theorem 2.7]. 

We treat the third case. From now on suppose that $\lambda\neq \lambda_i$ for any $i$. 
Pick $i$ so that $\lambda_i>\lambda_{i+1}$. Thus, $\Supp C_{i+1}$ does not 
contain any lc centre of $(X_{i+1}/Z,B_{i+1}+\lambda_{i+1}C_{i+1})$ because  $(X_{i+1}/Z,B_{i+1}+\lambda_{i}C_{i+1})$ is 
lc. Then, by replacing $(X/Z,B)$ with $(X_{i+1}/Z,B_{i+1})$ 
and $C$ with $\lambda_{i+1}C_{i+1}$ we may 
assume that no lc centre of $(X/Z,B+C)$ is inside $\Supp C$.  
Furthermore, using induction and the special termination (cf. [\ref{B}, Lemma 3.6]) we can assume that 
the log flips do not intersect $\rddown{B}$. 
Since in each step $K_{X_i}+B_i+\lambda C_i$ is anti-ample$/Z_i$, the sequence is also a sequence of 
$K_{X}+B+\lambda C$-flips.
By replacing $B$ with $B+\lambda C$, $C$ with $(1-\lambda)C$, and  $\lambda_i$ with 
$\frac{\lambda_i-\lambda}{1-\lambda}$, we may assume that $\lambda=0$.\\

\emph{Step 2.} By assumptions there is a log minimal model $(Y/Z,B_Y)$ for $(X/Z,B)$.
Let $\phi \colon X\bir Y/Z$ be the corresponding birational map. Since $K_{X_i}+B_i+\lambda_i C_i$ 
is nef$/Z$, we may add an ample$/Z$ $\R$-divisor $G^i$ so that $K_{X_i}+B_i+\lambda_i C_i+G^i$ 
becomes ample$/Z$, in particular, it is movable$/Z$. We can choose the $G^i$ so that 
$\lim_{i\to \infty} {G^i}_1=0$ in $N^1(X_1/Z)$ where ${G^i}_1$ is the birational transform of $G^i$ 
on $X_1=X$. Therefore, 
$$
K_{X}+B\equiv \lim_{i\to \infty} (K_{X_1}+B_1+\lambda_i C_1+G^i_1)/Z
$$ 
which implies that $K_{X}+B$ is a limit of movable$/Z$ $\R$-divisors.  

Let $f\colon W\to X$ and $g\colon W\to Y$ 
be a common log resolution of $(X/Z,B+C)$ and $(Y/Z,B_Y+C_Y)$ where $C_Y$ is the birational transform of $C$. 
By applying the negativity lemma to $f$, we see that 
$$
E:=f^*(K_X+B)-g^*(K_Y+B_Y)=\sum_D a(D,Y,B_Y)D-a(D,X,B)D
$$
is effective (cf. [\ref{B}, Remark 2.6]) where $D$ runs over the prime divisors on $W$. 
If $E\neq 0$, let $D$ be a component of $E$. If $D$ is not exceptional$/Y$, 
then it must be exceptional$/X$ otherwise $a(D,X,B)=a(D,Y,B_Y)$ and $D$ cannot be a component of $E$. 
By definition of log minimal models, $a(D,Y,B_Y)=0$ hence $a(D,X,B)=0$ which again shows that 
$D$ cannot be a component of $E$. Therefore, $E$ is exceptional$/Y$.\\  

\emph{Step 3.} 
Let $B_W$ be the birational transform of $B$ plus the reduced exceptional divisor of $f$, 
and let $C_W$ be the 
birational transform of $C$ on $W$. Pick a sufficiently small $\delta\ge 0$.
Take a general ample$/Z$ 
divisor $L$ so that $K_W+B_W+\delta C_W+L$ is dlt and nef$/Z$.  
Since $(X/Z,B)$ is lc, 
$$
E':=K_W+B_W-f^*(K_X+B)=\sum_D a(D,X,B)D\ge 0
$$
where $D$ runs over the prime exceptional$/X$ divisors on $W$. So,
$$
K_W+B_W+\delta C_W=f^*(K_X+B)+E'+\delta C_W=g^*(K_Y+B_Y)+E+E'+\delta C_W
$$ 
Moreover, $E'$ is also exceptional$/Y$ because for any prime divisor $D$ on $Y$ which is exceptional$/X$, 
$a(D,Y,B_Y)=a(D,X,B)=0$ hence $D$ cannot be a component of $E'$. 

On the other hand, since $Y$ is $\Q$-factorial, there are 
exceptional$/Y$ $\R$-divisors $F, F'$ on $W$ such that  $C_W+F\equiv 0/Y$ and $L+F'\equiv 0/Y$. 
Now run the LMMP$/Y$ on $K_W+B_W+\delta C_W$ with scaling of $L$ 
which is the same as the LMMP$/Y$ on $E+E'+\delta C_W$ with scaling of $L$. Let $\lambda_i'$ and 
$\lambda'=\lim_{i\to \infty} \lambda_i'$ be the 
corresponding numbers. 
If $\lambda'>0$, then by step 1 the LMMP terminates since $L$ is ample$/Z$. Since $W\to Y$ is birational, 
the LMMP terminates only when $\lambda_i'=0$ for some $i$ which implies that $\lambda'=0$, a contradiction.
Thus, $\lambda'=0$. On some model $V$ in the process of the LMMP, 
the pushdown of $K_W+B_W+\delta C_W+\lambda_i'L$, say 
\begin{equation*}
\begin{split}
K_V+B_V+\delta C_V+\lambda_i'L_V &\\ 
& \equiv E_V+E_V'+\delta C_V+\lambda_i'L_V\\ 
& \equiv  E_V+E_V'-\delta F_V-\lambda_i'F_V'/Y
\end{split}
\end{equation*}
is nef$/Y$. Applying the negativity lemma over $Y$ shows that $E_V+E_V'-\delta F_V-\lambda_i'F_V'\le 0$. But 
if $i\gg 0$, then $E_V+E_V'\le 0$  because 
$\lambda_i'$ and $\delta$ are sufficiently small. Therefore, $E_V=E_V'=0$ as $E$ and $E'$ are effective.\\

\emph{Step 4.} 
We prove that $\phi\colon X\bir Y$ does not contract any divisors. Assume otherwise and let 
$D$ be a prime divisor on $X$ contracted by $\phi$. Then $D^\sim$ the birational transform of 
$D$ on $W$ is a component of $E$ because by definition of log minimal models $a(D,X,B)<a(D,Y,B_Y)$.
Now, in step 3 take $\delta=0$. The LMMP contracts $D^\sim$  
since $D^\sim$ is a component of $E$ and $E$ is contracted. But this is not possible 
because  $K_{X}+B$ is a limit of movable$/Z$ $\R$-divisors and $D^\sim$ is not a component of 
$E'$ so the pushdown of $K_W+B_W=f^*(K_{X}+B)+E'$ cannot negatively intersect a general curve on $D^\sim/Y$. 
Thus $\phi$ does not contract divisors, in particular, any prime divisor on $W$ which is exceptional$/Y$ 
is also exceptional$/X$. 
Though $\phi$ does not contract divisors but $\phi^{-1}$ might contract divisors. 
The prime divisors contracted by $\phi^{-1}$ appear on $W$.\\

\emph{Step 5.} 
Now take $\delta>0$ in step 3 which is sufficiently small by assumptions.  By induction and the special termination, when we run the 
LMMP$/Y$ on $K_W+B_W+\delta C_W$ with scaling of $L$, the extremal rays contracted in the process 
do not intersect $\rddown{B_W}$, after finitely many steps. On the other hand, 
since $\phi$ does not contract divisors, every exceptional$/Y$ prime divisor
on $W$ is a component of $\rddown{B_W}$. Therefore, the LMMP terminates because 
it is an LMMP on the exceptional$/Y$ $\R$-divisor $E+E'-\delta F$.
So, we get a model $Y'$ on which 
the pushdown of $K_W+B_W+\delta C_W$, say $K_{Y'}+B_{Y'}+\delta C_{Y'}$, is nef$/Y$. 
By step 3, $K_{Y'}+B_{Y'}\equiv E_{Y'}+E'_{Y'}=0/Y$ where $E_{Y'}$ and $E'_{Y'}$ are the birational 
transforms of $E$ and $E'$ on $Y'$, respectively.
Therefore, $(Y'/Z,B_{Y'})$ is a dlt crepant model of $(Y/Z,B_Y)$.\\

\emph{Step 6.} As in step 3, 
$$
E'':=K_W+B_W+C_W-f^*(K_X+B+C)=\sum_D a(D,X,B+C)D\ge 0
$$
is exceptional$/X$ where $D$ runs over the prime exceptional$/X$ divisors on $W$.
 So, by induction and the special termination, the LMMP$/X$ on 
$K_W+B_W+C_W\equiv E''/X$ with scaling of suitable ample$/Z$ divisors terminate because 
every component of $E''$ is also a component of $\rddown{B_W}$.
So, we get a crepant dlt model $(X'/Z,B'+C')$ of $(X/Z,B+C)$ 
where $K_{X'}+B'$ is the pullback of $K_{X}+B$ and $C'$ is the pullback of $C$. 
In fact, $X'$ and $X$ are isomorphic outside the lc centres of $(X/Z,B+C)$ because the 
prime exceptional$/X$ divisors on $X'$ are exactly the pushdown of the prime 
exceptional$/X$ divisors $D$ on $W$ with $a(D,X,B+C)=0$, that is, those which are not 
components of $E''$. Since 
$\Supp C$ does not contain any lc centre of $(X/Z,B+C)$ by step 1, 
$(X'/Z,B')$ is a crepant dlt model of $(X/Z,B)$ and $C'$ is just
the birational transform of $C$. 
Note that the prime exceptional divisors of  $\phi^{-1}$ 
are not contracted$/X'$ since their log discrepancy with respect to $(X/Z,B)$ are all $0$, and so their birational 
transforms are not components of $E''$.\\

\emph{Step 7.} Remember that $X_1=X$, $B_1=B$, and $C_1=C$. Similarly, put $X_1':=X'$, $B_1':=B'$, and $C_1':=C'$. 
Since $K_{X_1}+B_1+\lambda_1C_1\equiv 0/Z_1$, 
$K_{X_1'}+B_1'+\lambda_1C_1'\equiv 0/Z_1$. Run the LMMP$/Z_1$ on $K_{X_1'}+B_1'$ 
with scaling of $\lambda_1 C_1'$. Since the exceptional locus of $X_1\to Z_1$ 
does not intersect any lc centre of $(X_1/Z,B_1)$ by step 1, and since $X_1'$ and $X_1$ are isomorphic 
outside the lc centres of $(X_1/Z,B_1)$, the LMMP consists of just one log flip $X_1'\bir X_2'/Z_1'$ which is the lifting 
of the log flip $X_1\bir X_2/Z_1$. Moreover, $(X_2'/Z,B_2')$ is a crepant dlt model of $(X_2/Z,B_2)$ where 
$B_2'$ is the birational transform of $B_1'$. 
We can continue this process to lift the original sequence to a sequence $X_i'\bir X_{i+1}'/Z_i'$. 

Note that $Y'\bir X'$ does not 
contract divisors: if $D$ is a prime divisor on $Y'$ which is exceptional$/X'$, then it is exceptional$/X$ and 
so it is exceptional$/Y$ by step 6; but then $a(D,Y,B_Y)=0=a(D,X,B)$ and again by step 6 such divisors 
are not contracted$/X'$, a cotradiction. Thus, $(Y'/Z,B_{Y'})$ of step 5 is a log birational model of $(X'/Z,B')$ 
because $B_{Y'}$ is the birational transform of $B'$.
On the other hand, assume that $D$ is a 
prime divisor on $X'$ which is exceptional$/Y'$. Since $X\bir Y$ does not contract divisors by step 4, $D$ is exceptional$/X$. 
In particular, $a(D,X',B')=a(D,X,B)=0$; in this case 
$a(D,Y,B_Y)=a(D,Y',B_{Y'})>0$ 
otherwise $D$ could not be contracted$/Y'$ by the LMMP of step 5 which started on $W$ because 
the birational transform of $D$ would not be a component of $E+E'+\delta C_W$. 
So, $(Y'/Z,B_{Y'})$ is actually a log minimal model of $(X'/Z,B')$.
Therefore, as in step 4, $X'\bir Y'$ does not contract divisors which implies that $X'$ and $Y'$ are 
isomorphic in codimension one. Now replace the old sequence $X_i\bir X_{i+1}/Z_i$
with the new one $X_i'\bir X_{i+1}'/Z_i'$ and replace $(Y/Z,B_Y)$ with $(Y'/Z,B_{Y'})$.
So, from now on we can assume that $X, X_i$ and $Y$ 
are all isomorphic in codimension one. In addition, by step 5, we can also assume that 
$(Y/Z,B_{Y}+\delta C_{Y})$ is dlt for some $\delta>0$.\\

\emph{Step 8.}
Let $A\ge 0$ be a reduced divisor on $W$ whose components 
are general ample$/Z$ divisors such that they generate $N^1(W/Z)$. By step 6, $(X_1/Z,B_1+C_1)$ is obtained by running 
a specific LMMP on 
$K_W+B_W+C_W$. Every step of this LMMP is also a step of an LMMP on $K_W+B_W+C_W+\ep A$ 
for any sufficiently small $\ep>0$, in particular, $(X_1/Z,B_1+C_1+\ep A_1)$ is dlt where $A_1$ is the birational transform of $A$. 
For similar reasons, we can choose $\ep$ so that 
$(Y/Z,B_Y+\delta C_Y+\ep A_Y)$ is also dlt. On the other hand, by Proposition \ref{p-rays}, 
perhaps after replacing $\delta$ and $\ep$ with 
smaller positive numbers, we may assume that if $0\le \delta'\le \delta$ and $0\le A_Y'\le A_Y$, then any 
 LMMP$/Z$ on $K_Y+B_Y+\delta' C_Y+A_Y'$, consists of only a sequence of log flips which are flops with respect to $(Y/Z,B_Y)$. 
Note that since $K_Y+B_Y+\delta' C_Y+A_Y'$ is a limit of movable$/Z$ $\R$-divisors, no divisor is contracted 
by such an LMMP.\\

\emph{Step 9.} Fix some $i\gg 0$ so that $\lambda_i<\delta$. Then, by Proposition \ref{p-rays}, there is $0<\tau \ll \ep$ 
such that $(X_i/Z,B_i+\lambda_i C_i+\tau A_i)$ is dlt and such that if we run the LMMP$/Z$ 
on  $K_{X_i}+B_i+\lambda_i C_i+\tau A_i$ with scaling of some ample$/Z$ divisor, then it will be a 
sequence of log flips which would be a sequence of flops with respect to $(X_i/Z,B_i+\lambda_i C_i)$. 
Moreover, since the components of $A_i$ generate $N^1(X_i/Z)$, we can assume that there is an ample$/Z$ $\R$-divisor 
$H\ge 0$ such that $\tau A \equiv H+H'/Z$ where $H'\ge 0$ and $(X_i/Z,B_i+\lambda_i C_i+H+H')$
is dlt. Hence the LMMP terminates by step 1
and we get a model $T$ on which both $K_T+B_T+\lambda_i C_T$ 
and $K_T+B_T+\lambda_i C_T+\tau A_T$ are nef$/Z$.  Again since the components of $A_T$ generate $N^1(T/Z)$,  
there is $0\le A'_T\le \tau A_T$ so that $K_T+B_T+\lambda_i C_T+ A_T'$ is ample$/Z$ and $\Supp A_T'=\Supp A_T$.
Now run the LMMP$/Z$ on 
$K_Y+B_Y+\lambda_i C_Y+A'_Y$ with scaling of some ample$/Z$ divisor  where $A_Y'$ 
is the birational tranform of $A_T'$. The LMMP terminates for reasons similar to the above and we end up
with $T$ since $K_T+B_T+\lambda_i C_T+ A_T'$ is ample$/Z$. Moreover, the LMMP consists of only 
log flips which are flops with respect to $(Y/Z,B_Y)$ by Proposition \ref{p-rays} hence $K_T+B_T$ 
will also be nef$/Z$. So, by replacing $Y$ with $T$ we could assume that $K_Y+B_Y+\lambda_i C_Y$ 
is nef$/Z$. In particular, $K_Y+B_Y+\lambda_j C_Y$ is nef$/Z$ for any $j\ge i$ since $\lambda_j\le \lambda_i$.\\ 

\emph{Step 10.} Pick $j> i$ so that $\lambda_j<\lambda_{j-1}\le \lambda_i$ and let $r\colon U\to X_j$ and $s\colon U\to Y$ 
be a common resolution. Then, we have 

\begin{equation*}
\begin{split}
r^*(K_{X_j}+B_j+\lambda_j C_j) & = s^*(K_Y+B_Y+\lambda_j C_Y)\\ 
r^*(K_{X_j}+B_j) & \gneq s^*(K_Y+B_Y)\\
r^*C_j & \lneq  s^* C_Y
\end{split}
\end{equation*}
where the first equality holds because both $K_{X_j}+B_j+\lambda_j C_j$ and $K_Y+B_Y+\lambda_j C_Y$ 
are nef$/Z$ and $X_j$ and $Y$ are isomorphic in codimension one, the second inequality holds 
because $K_Y+B_Y$ is nef$/Z$ but $K_{X_j}+B_j$ is not nef$/Z$, and the third follows from the 
other two. Now

\begin{equation*}
\begin{split}
r^*(K_{X_j}+B_j+\lambda_{j-1} C_j) & \\
 & = r^*(K_{X_j}+B_j+\lambda_{j} C_j)+ r^*(\lambda_{j-1}-\lambda_j)C_j\\
& \lneq  s^*(K_Y+B_Y+\lambda_{j} C_Y)+s^*(\lambda_{j-1}-\lambda_j)C_Y\\
& =  s^*(K_Y+B_Y+\lambda_{j-1} C_Y)
\end{split}
\end{equation*}

However, since $K_{X_j}+B_j+\lambda_{j-1} C_j$ 
and  $K_Y+B_Y+\lambda_{j-1} C_Y$ are both nef$/Z$, we have 
$$
r^*(K_{X_j}+B_j+\lambda_{j-1} C_j) =s^*(K_Y+B_Y+\lambda_{j-1} C_Y)
$$ 
This is a contradiction and the sequence of log flips terminates as claimed.
$\Box$\\

{\textbf{Proof of Corollary \ref{c-1}.}}
Let $H\ge 0$ be an ample$/Z$ divisor such that $K_X+B+H$ is dlt and ample$/Z$. Now run the 
LMMP$/Z$ on $K_X+B$ with scaling of $H$. By Theorem \ref{t-main2}, the LMMP terminates 
with a log minimal model or a Mori fibre space $(Y/Z,B_Y)$. The claim that $Y\bir X$ does not 
contract divisors is obvious. 
$\Box$\\

\begin{lem}\label{l-special}
Assume the minimal model conjecture (\ref{mmodel}) in dimension $d$ for pseudo-effective 
$\Q$-factorial dlt pairs. Let $(X/Z,B+C)$ be a 
$\Q$-factorial lc pair of dimension $d+1$ such that 
\begin{enumerate}
\item $K_X+B+C$ is nef$/Z$,
\item $B,C\ge 0$,
\item $(X/Z,B)$ is dlt,
\item $K_X+B\equiv_Z M\ge 0$ where $\alpha M=M'+C$ for some $\alpha>0$ and $M'\ge 0$ supported in $\Supp \rddown{B}$.
\end{enumerate} 

Then, we can run an LMMP$/Z$ on $K_X+B+C$ with scaling of $C$ 
which terminates.
\end{lem}
{\textbf{Proof.}}
By Theorem \ref{t-main2}, [\ref{HM}, Assumption 5.2.3] is satisfied in dimension $d$ which implies that pl flips exist in dimension $d+1$ 
by the main result of [\ref{HM}] (cf. [\ref{BP}, Theorem 2.9]). Alternatively, we can simply borrow 
the existence of log flips from  
[\ref{BCHM}][\ref{BP}]. So, in any case we can run the LMMP$/Z$ on $K_X+B$ with scaling of $C$ by [\ref{B}, Lemma 3.1] 
because we only need pl flips. We may 
assume that any LMMP$/Z$ on $K_X+B$ with scaling of $C$ consists of only log flips.   

If $M'=0$, then $K_X+B+C\equiv \frac{1}{\alpha}C+C/Z$ which implies that 
$C$  and $K_X+B$ are nef$/Z$ hence we are done.  So, from now on we assume 
that $M'\neq 0$. 

By the assumptions, $\Supp M\subseteq \Supp (B+C)$ hence there is a sufficiently small $\tau>0$ such that 
$$
\Supp (B+C-\tau M-\tau C)=\Supp (B+C)
$$
Put $B'=B-\frac{\tau}{\alpha} M'$ and $C'=C-\tau (\frac{1}{\alpha}+1)C$ so that 
$$
K_X+B'+C'\equiv M+C-\frac{\tau}{\alpha} M'-(\frac{\tau}{\alpha}+\tau)C=M+C-\tau (M+C)/Z
$$
In particular, $K_X+B'+C'$ is nef$/Z$. Let $\delta$ be as in  Proposition \ref{p-rays} 
chosen for the pair $(X/Z,B'+C')$ where we take $V$ to be the space $V=\{rM' \mid r\in \R\}$. 
Take $a>0$ so that $a\alpha \ll \tau$, $||aM'||<\delta$, 
$B'':=B-a M'\ge 0$ has the same support as $B$, and $C''=C-(a+a\alpha)C\ge 0$ 
has the same support as $C$. Now 
$$
K_X+B''+C''\equiv M+C-a M'-(a+a\alpha)C=M+C- a\alpha(M+C)/Z
$$
and $\Supp M'\subseteq \Supp \rddown{B=B''+a M'}$. In particular, $K_X+B''+C''$ is nef$/Z$. 

Let $H\ge 0$ be an ample$/Z$ divisor such that 
$K_X+B+C''+H$ 
is dlt and ample$/Z$. Now run the LMMP$/Z$ on $K_X+B+C''$ with scaling of $H$ 
and assume that we get a sequence $X_i\bir X_{i+1}$ of log flips and divisorial contractions 
corresponding to extremal rays $R_i$. For each $i$, we have 
$$
0>(K_{X_i}+B_i+C_i'')\cdot R_i=(1-a\alpha)(M_i+C_i)\cdot R_i+aM'_i\cdot R_i
$$ 
where as usual the subscript $i$ for divisors stands for birational transform on $X_i$. 
By induction on $i$, we may assume that $K_{X_i}+B_i''+C_i''$ is nef$/Z$ which also means that 
$K_{X_i}+B_i'+C_i'$ is nef$/Z$. So $M'_i\cdot R_i<0$ and 
\begin{equation*}
\begin{split}
(K_{X_i}+B_i'+C_i'+a M'_i)\cdot R_i &=(1-\tau)(M_i+C_i)\cdot R_i+aM'_i\cdot R_i\\
&<(1-a\alpha)(M_i+C_i)\cdot R_i+aM'_i\cdot R_i<0
\end{split}
\end{equation*}
which implies that $(K_{X_i}+B_i'+C_i')\cdot R_i=0$, by construction, and in turn we get $(K_{X_i}+B_i''+C_i'')\cdot R_i=0$. Thus, 
$C_i\cdot R_i>0$ and $(K_{X_i}+B_i+C_i)\cdot R_i=0$. 
So, the above LMMP is an LMMP$/Z$ on $K_X+B$ with scaling of $C$. Since $H$ is ample$/Z$, 
the LMMP terminates by the special termination and Theorem \ref{t-main2} because the LMMP is a 
$(-M')$-LMMP and  $\Supp M'\subseteq \Supp \rddown{B}$. Thus, for some $i$, 
 $K_{X_i}+B_i+C_i''=K_{X_i}+B_i+(1-a-a\alpha)C_i$ is nef$/Z$. 

Now replace $(X/Z,B)$ 
with $(X_i/Z,B_i)$, $C$ with $C_i''=(1-a-a\alpha)C_i$, $M$ with $M_i$, $M'$ with $(1-a-a\alpha)M'_i$, $\alpha$ 
with $\alpha(1-a-a\alpha)$, and continue the process by starting from the beginning. This process stops again 
by the special termination and Theorem \ref{t-main2}.

The underlying idea is that there is an LMMP$/Z$ on $K_X+B$ with scaling of $C$ such that the 
corresponding numbers $\lambda_i$ and $\lambda$ satisfy the property $\lambda\neq \lambda_i$ for any $i$ 
and this allows us to use the special termination and apply Theorem \ref{t-main2} in lower dimension.
$\Box$\\

{\textbf{Proof of Corollary \ref{c-2}.}}
Let $(X/Z,B)$ be an effective lc pair of dimension $d+1$. 
By [\ref{B}, Proposition 3.4], existence of pl flips in dimension $d+1$ and 
the special termination with scaling in dimension $d+1$ for $\Q$-factorial dlt pairs 
implies the existence of a log minimal model 
for $(X/Z,B)$. As mentioned in the proof of Lemma \ref{l-special}, existence of pl flips in dimension $d+1$ 
follows from the assumptions. However, we have not derived termination with scaling in dimension $d$ from our 
assumptions when $\lambda=\lambda_i$ for some $i$. But this is not a problem since we 
can use Lemma \ref{l-special}. We analyse the various places in the proof of [\ref{B}, Proposition 3.4] 
where the special termination is needed. 

In step $1$ of the proof of  [\ref{B}, Proposition 3.4]  we need to have special termination with scaling 
of an ample$/Z$ $\R$-divisor for a certain sequence of log flips. This follows from our assumptions by  
Theorem \ref{t-main2}. In steps $3,4$, and $5$ 
we need the special termination for some LMMP with scaling in a situation as follows: $(X/Z,B+C)$ is log 
smooth, $B,C\ge 0$, $K_X+B\equiv_Z M\ge 0$, $\alpha M=M'+C$ for some $\alpha>0$, 
$M'\ge 0$ is supported in $\Supp \rddown{B}$, and 
$(Y/Z,B_Y+C_Y)$ is a log minimal model of  $(X/Z,B+C)$ where $B_Y$ is the birational 
transform of $B$ plus the reduced exceptional divisor of $Y\bir X$ and $C_Y$ is just the birational 
transform of $C$. Here we want to run an LMMP$/Z$ on $K_Y+B_Y$ with scaling of $C_Y$ which terminates.
Let $f\colon W\to X$ and $g\colon W\to Y$ be a common log resolution. By the arguments in step 
2 of the proof of Theorem \ref{t-main2}, we can write $f^*(K_X+B+C)=g^*(K_Y+B_Y+C_Y)+E$ 
where $E$ is effective, and exceptional$/Y$. So, 
$$
f^*(M+C)=f^*(\frac{1}{\alpha}M'+\frac{1}{\alpha}C+C) \equiv_Z g^*(K_Y+B_Y+C_Y)+E
$$ 
and 
$$
g_*f^*(\frac{1}{\alpha}M'+\frac{1}{\alpha}C+C)\equiv_Z K_Y+B_Y+C_Y
$$
Now put $M_Y:=g_*f^*(\frac{1}{\alpha}M'+\frac{1}{\alpha}C+C)-C_Y$ and $M_Y':=\alpha M_Y-C_Y$ 
so that $K_Y+B_Y\equiv M_Y/Z$ and  $\alpha M_Y=M_Y'+C_Y$. By construction, every component of $M_Y'$ is either the birational 
trasnform of a component of $M'$ or it is an exceptional divisor of $Y\bir X$ which in any case would be a 
component of $\rddown{B_Y}$. Now simply apply Lemma \ref{l-special} to the data: $(Y/Z,B_Y+C_Y)$, 
$M_Y$, $\alpha$, and $M_Y'$. 

In step $6$  of the proof of  [\ref{B}, Proposition 3.4]  we need 
special termination to be able to apply [\ref{B}, Lemma 3.3]. However, the proof of [\ref{B}, Lemma 3.3] 
only needs the special termination with scaling 
of an ample$/Z$ $\R$-divisor applied to a certain sequence of log flips which again follows from our assumptions by Theorem \ref{t-main2}. 
$\Box$\\

{\textbf{Proof of Theorem \ref{t-main}.}}
We use induction on $d$ so assume that the theorem holds in dimension $d-1$. In particular, 
we may assume that the minimal model conjecture (\ref{mmodel}) holds in dimension $d-1$.
Let $(X/Z,B)$ be a lc pair of dimension $d$. We may assume that $(X/Z,B)$ is $\Q$-factorial 
dlt by replacing it with a $\Q$-factorial dlt crepant model. To construct such a model (cf. step 
6 of the proof of Theorem \ref{t-main2}) we only need 
 the special termination with scaling 
of an ample$/Z$ $\R$-divisor applied to a certain sequence of log flips which follows from the minimal model 
conjecture in dimension $d-1$ and  Theorem \ref{t-main2}.
If $K_X+B$ is not pseudo-effective$/Z$, then by [\ref{BCHM}][\ref{BP}] there is a Mori fibre space 
for $(X/Z,B)$. If  $K_X+B$ is pseudo-effective$/Z$, then by Conjecture \ref{c-wa}, it is effective, that is, 
there is $M\ge 0$ such that $K_X+B\equiv M/Z$. Now the result follows from Corollary \ref{c-2}. 

The statement concerning  $\Q$-factorial dlt $(X/Z,B)$ follows from Corollary \ref{c-1}, that is, we can 
run the LMMP$/Z$ on $K_X+B$ with 
scaling of some ample$/Z$ $\R$-divisor which will end up with a log minimal model or a Mori fibre space. 
$\Box$\\


\vspace{2cm}

\flushleft{DPMMS}, Centre for Mathematical Sciences,\\
Cambridge University,\\
Wilberforce Road,\\
Cambridge, CB3 0WB,\\
UK\\
email: c.birkar@dpmms.cam.ac.uk


\begin{thebibliography}{99}

\bibitem{}\label{B}  {C. Birkar; {\emph{On existence of log minimal models.}} arXiv:0706.1792v3. }

\bibitem{}\label{B2}  {C. Birkar; {\emph{Log minimal models according to Shokurov.}} arXiv:0804.3577v1. }

\bibitem{}\label{BCHM}  {C. Birkar, P. Cascini, C. Hacon, J. M$^c$Kernan; {\emph{Existence of minimal models
for varieties of log general type.}}  arXiv:math/0610203v2.}

\bibitem{}\label{BP}  {C. Birkar, M. P\u{a}un; {\emph{Minimal models, flips and finite generation : a tribute to 
V.V. SHOKUROV and Y.-T. SIU.}}  arXiv:0904.2936v1.}

\bibitem{}\label{HM}  {C. Hacon, J. M$^c$Kernan; {\emph{  Extension theorems and the existence of flips. In}} 
Flips for $3$-folds and $4$-folds, Oxford University Press (2007).}

\bibitem{}\label{Kawamata}  {Y. Kawamata; {\emph{Flops connect minimal models.}}
Publ. RIMS, Kyoto Univ. 44 (2008), 419-423.}

\bibitem{}\label{log-flips}  {V.V. Shokurov; {\emph{3-fold log flips.}}
With an appendix in English by Yujiro Kawamata.
Russian  Acad. Sci. Izv. Math.  40  (1993),  no. 1, 95--202.}

\bibitem{}\label{antican}  {V.V. Shokurov; \emph{ Anticanonical boundedness for curves.} Appendix to
V.V. Nikulin "Hyperbolic reflection group methods and algebraic varieties"
in Higher dimensional complex varieties (Trento, June 1994) ed.
Andreatta M., and Peternell T. Berlin: New York: de Gruyter (1996), 321-328.}

\bibitem{}\label{ordered}  {V.V. Shokurov; {\emph{Letters of a bi-rationalist VII: Ordered termination.}}
arXiv:math/0607822v2.}

\end{thebibliography}
\end{document}